# A space efficient flexible pivot selection approach to evaluate determinant and inverse of a matrix

Hafsa Athar Jafree[1], Muhammad Imtiaz[2], Syed Inayatullah[3],
Fozia Hanif Khan[4] and Tajuddin Nizami[5]

## ABSTRACT

This paper presents new approaches for finding the determinant and inverse of a matrix. The choice of pivot selection is kept arbitrary and can be made according to the users need. So the ill conditioned matrices can be handled easily. The algorithms are more efficient as they save unnecessary data storage by reducing the order of the matrix after each iteration in the computation of determinant and incorporating dictionary notation (Chvatal, 1983) in the computation of inverse matrix. These algorithms are highly class room oriented and unlike the matrix inversion method (Khan, Shah, & Ahmad, 2010) the presented algorithm doesn't need any kind of permutations or inverse permutations.

[1] Department of Mathematical Sciences, University of Karachi, Karachi,Pakistan. hajafree@uok.edu.pk
[2] Department of Mathematical Sciences, University of Karachi, Karachi, Pakistan. imtiaz@uok.edu.pk
[3] Department of Mathematical Sciences, University of Karachi, Karachi, Pakistan. inayat@uok.edu.pk
[4] Department of Mathematics, Sir Syed University of Engineering and Technology, Karachi, Pakistan. nizamitaj@yahoo.com
[5] Department of Mathematics, Iqra University, Karachi, Pakistan. nizamitaj@yahoo.com

# 1. INTRODUCTION

Determinant is associated with the system of linear equations. It provides us with a foresight about the nature of solution of a given system of linear equations.

Maclaurin (1748) published the first result on determinant of $2 \times 2$ and $3 \times 3$ systems, which was later generalized for $n \times n$ systems by Cramer (1750). Laplace used the word 'resultant' for determinant and proposed a method now known as Laplace expansion. Gauss (1801) used the word determinant for discussing the properties of quadratic form. Cauchy (1815) used the word determinant for the first time in the modern context. In 1866, Dodgson (Rice & Torrence, 2006) presented the method for finding the determinant of $n \times n$ systems which he named as the "Method of condensation".

The existence of matrix inverse also depends on its determinant. (Ahmad & Khan, 2010), (Khan, Shah, & Ahmad, 2010) proposed algorithms for the calculation of matrix inverse which is a streamlined form of the Gaussian method. In section 2 of this paper, we have presented a new algorithm to find determinant and in section 4 we have also presented another easy to manage way for evaluating the inverse of a matrix, which doesn't require any kind of permutations and inverse permutations.

# 2. METHOD FOR FINDING DETERMINANT

**Problem**

Find determinant of the matrix $\mathbf{A} = [a_{ij}]$ of order $n \times n$.

**Algorithm**

**Step 1:** Set $d:=1$,

Step 2: set $P:=\{1,2,3,\ldots,n\}$

**Step 3:** For any $p \in P$, $L:=\{k: a_{pk} \neq 0, k \in P\}$.

**Step 4:** If $L=\emptyset$ then $d=0$ and go to step 6

$$\text{Otherwise } a_{ij} = a_{ij} - \frac{a_{ik} \times a_{pj}}{a_{pk}} \qquad \forall i \in P - \{p\} \text{ and } j \in P - \{k\}$$

**Step 4:** $d := (-1)^{p+k} d$

**Step 5:** Reduce the order of A by removing $p^{th}$ row and $k^{th}$ column of A, $n := n - 1$. If $n \neq 0$ go to step 1.

**Step 6:** $det(A) = d$. Exit.

**Example**

Consider the matrix $\begin{bmatrix} 2 & 5 & 3 & 2 \\ 4 & 10 & 1 & 7 \\ 1 & 5 & 2 & 2 \\ 2 & 1 & 2 & 1 \end{bmatrix}$

**Iteration 1:**

$P := \{1,2,3,4\}$, here we take $p = 1$ so, $L := \{1,2,3,4\}$. Taking $k = 1$

$A := \begin{bmatrix} 0 & -5 & 3 \\ 5/2 & 1/2 & 0 \\ -4 & -1 & -1 \end{bmatrix}$

$d := (-1)^{1+1} \times 2 = 2$

**Iteration 2:**

$P := \{1,2,3\}$, here we take $p = 2$ so, $L := \{2,3\}$. Taking $k = 2$

$A := \begin{bmatrix} 5/2 & 3/10 \\ -4 & -8/5 \end{bmatrix}$

$d := 2 \times (-1)^{1+2} \times (-5) = 10$

**Iteration 3:**

$P := \{1,2\}$, here we take $p = 1$ so, $L := \{1,2\}$. Taking $k = 1$

$A := [-28/25]$

$d := 10 \times (-1)^{1+1} \times \frac{5}{2} = 25$

**Iteration 4:**

$P := \{1\}$, here we take $p = 1$ so, $L := \{1\}$. Taking $k = 1$

$d := 25 \times (-1)^{1+1} \times \left(-\frac{28}{25}\right) = -28$

Hence determinant of given matrix is $-28$

## 3. INVERSE MATRIX

Consider a matrix of order $n$. To evaluate the inverse the matrix $\mathbf{A} = [a_{ij}]$ one must solve the following $n$ system of equations, for $x_1, x_2, \ldots, x_n$.

$$\begin{matrix} a_{11}x_1 & a_{12}x_2 & \cdots & a_{1n}x_n = & y_1 \\ a_{21}x_1 & a_{22}x_2 & \cdots & a_{2n}x_n = & 0 \\ \vdots & \vdots & & \vdots \vdots & \vdots \\ a_{n1}x_1 & a_{n2}x_2 & \cdots & a_{nn}x_n = & 0 \end{matrix} ,\quad \begin{matrix} a_{11}x_1 & a_{12}x_2 & \cdots & a_{1n}x_n = & 0 \\ a_{21}x_1 & a_{22}x_2 & \cdots & a_{2n}x_n = & y_2 \\ \vdots & \vdots & & \vdots \vdots & \vdots \\ a_{n1}x_1 & a_{n2}x_2 & \cdots & a_{nn}x_n = & 0 \end{matrix} ,\ldots,$$

$$\begin{matrix} a_{11}x_1 & a_{12}x_2 & \cdots & a_{1n}x_n = & 0 \\ a_{21}x_1 & a_{22}x_2 & \cdots & a_{2n}x_n = & 0 \\ \vdots & \vdots & & \vdots \vdots & \vdots \\ a_{n1}x_1 & a_{n2}x_2 & \cdots & a_{nn}x_n = & y_n \end{matrix}$$

These matrices can also be written compactly in the following augmented matrix form as,

$$\begin{matrix} x_1 & x_2 & \cdots & x_n & y_1 & y_2 & \cdots & y_n \\ \begin{bmatrix} a_{11} & a_{12} & \cdots & a_{1n} & 1 & 0 & \cdots & 0 \\ a_{21} & a_{22} & \cdots & a_{2n} & 0 & 1 & \cdots & 0 \\ \vdots & \vdots & \vdots & \vdots & \vdots & \vdots & \vdots & \vdots \\ a_{n1} & a_{n2} & \cdots & a_{nn} & 0 & 0 & \cdots & 1 \end{bmatrix} \end{matrix}$$

Gaussian elimination is could be applied to obtain the inverse of the matrix. Here, in our technique we have assumed $B = \{y_1, y_1, \cdots, y_n\}$ as the basis and $N = \{x_1, x_2, \cdots, x_n\}$ as the non-basis of the current matrix. The objective is to convert the basic variables into non-basic variables and vice versa. For this purpose we may follow the dictionary notation developed by (Chvatal, 1983). After removal of basic columns from the above matrix we get the following dictionary form with basis $B$ and non-basis $N$:

$$\begin{array}{c c c c c} & x_1 & x_2 & \cdots & x_n \\ y_1 & \begin{bmatrix} a_{11} & a_{12} & \cdots & a_{nn} \\ a_{21} & a_{22} & \cdots & a_{2n} \\ \vdots & \vdots & \vdots & \vdots \\ a_{n1} & a_{n2} & \cdots & a_{nn} \end{bmatrix} \\ y_2 & & & & \\ \vdots & & & & \\ y_n & & & & \end{array}$$

The following pivot operation may be applied to enter $x_i$ into basis $B$ and $y_j$ into non-basis $N$,

i) Divide the pivot row by pivot element and the pivot column by negative of the pivot element (except the pivot element).

ii) The remaining $(n-1)^2$ elements are determined by the formula as mentioned in the following algorithm.

iii) Reciprocate the pivot element.

If the number of pivot elements is equal to the order of the matrix the resulting matrix gives the inverse otherwise we may conclude that the inverse does not exist.

## 4. METHOD FOR FINDING INVERSE

**Problem**

Find the inverse of matrix $\mathbf{A}=[a_{ij}]$ of order $n \times n$.

**Algorithm**

**Step 1:** Set $H:=\{1,2,3,\ldots,n\}$, $B:=\{y_h\}$ and $N:=\{x_h\}$, $\forall\, h \in H$. Construct dictionary of the matrix $\mathbf{A}$, i.e. D($\mathbf{A}$).

**Step 2:** Set $P:=\{p: y_p \in B\}$

**Step 3:** If $P = \emptyset$, go to step 6.

Otherwise, $L:=\{k: a_{pk} \neq 0,\ x_k \in N\}$

**Step 4:** If $L=\emptyset$ then inverse does not exist. Exit

Otherwise, for any $p \in P$, $k \in L$

$$m_i := -\frac{a_{ik}}{a_{pk}} \qquad \forall i \in H - \{p\}$$

$$a_{ij} := a_{ij} + a_{pj} \times m_i \qquad \forall i \in H - \{p\},\ \forall j \in H - \{k\}$$

$$a_{pj} := \frac{a_{pj}}{a_{pk}} \qquad \forall j \in H - \{k\}$$

$$a_{ik} := m_i \qquad \forall i \in H - \{p\}$$

$$a_{pk} := \frac{1}{a_{pk}}$$

**Step 5:** $B := B + \{x_k\} - \{y_p\}$ and $N := N - \{x_k\} + \{y_p\}$. Go to step 2.

**Step 6:** $Inv(\mathbf{A}) = [a_{ij}, \forall\, i, j \in H: x_i \in B, y_j \in N]$. Exit.

**Example**

Consider $\mathbf{A} = \begin{bmatrix} 2 & 5 & 3 & 2 \\ 4 & 10 & 1 & 7 \\ 1 & 5 & 2 & 2 \\ 2 & 1 & 2 & 1 \end{bmatrix}$

$$D(\mathbf{A}) := \begin{array}{c} \\ y_1 \\ y_2 \\ y_3 \\ y_4 \end{array} \begin{array}{cccc} x_1 & x_2 & x_3 & x_4 \\ \begin{bmatrix} 2 & 5 & 3 & 2 \\ 4 & 10 & 1 & 7 \\ 1 & 5 & 2 & 2 \\ 2 & 1 & 2 & 1 \end{bmatrix} \end{array}$$

Here $N := \{x_1, x_2, x_3, x_4\}$, $B := \{y_1, y_2, y_3, y_4\}$

**Iteration 1:**

$H := \{1,2,3,4\}$, $P := \{1,2,3,4\}$, taking $p=1$ we get $L := \{1,2,3,4\}$. Taking $k=1$

$$D(\mathbf{A}) := \begin{array}{c} \\ x_1 \\ y_2 \\ y_3 \\ y_4 \end{array} \begin{array}{cccc} y_1 & x_2 & x_3 & x_4 \\ \begin{bmatrix} 1/2 & 5/2 & 3/2 & 1 \\ -2 & 0 & -5 & 3 \\ -1/2 & 5/2 & 1/2 & 0 \\ -1 & -4 & -1 & -1 \end{bmatrix} \end{array}$$

$N := \{y_1, x_2, x_3, x_4\}$, $B := \{x_1, y_2, y_3, y_4\}$

**Iteration 2:**

$P := \{2,3,4\}$, taking $p=2$ we get $L := \{3,4\}$. Taking $k=3$

$$D(\mathbf{A}) := \begin{array}{c} \\ x_1 \\ x_3 \\ y_3 \\ y_4 \end{array} \begin{array}{c} y_1 \quad\; x_2 \quad\; y_2 \quad\; x_4 \end{array} \\ \begin{bmatrix} -1/10 & 5/2 & 3/10 & 19/10 \\ 2/5 & 0 & -1/5 & -3/5 \\ -7/10 & 5/2 & 1/10 & 3/10 \\ -3/5 & -4 & -1/5 & -8/5 \end{bmatrix}$$

$N := \{y_1, x_2, y_2, x_4\}$, $B := \{x_1, x_3, y_3, y_4\}$

**Iteration 3:**

$P := \{3,4\}$, taking $p=3$ we get $L := \{2,4\}$. Taking $k=2$

$$D(\mathbf{A}) := \begin{array}{c} \\ x_1 \\ x_3 \\ x_2 \\ y_4 \end{array} \begin{array}{c} y_1 \quad\;\; y_3 \quad\;\; y_2 \quad\;\; x_4 \end{array} \\ \begin{bmatrix} 3/5 & -1 & 1/5 & 8/5 \\ 2/5 & 0 & -1/5 & -3/5 \\ -7/25 & 2/5 & 1/25 & 3/25 \\ -43/25 & 8/5 & -1/25 & -28/25 \end{bmatrix}$$

$B := \{x_1, x_3, x_2, y_4\}$, $N := \{y_1, y_3, y_2, x_4\}$

**Iteration 4:**

$P := \{4\}$, taking $p=4$ we get $L := \{4\}$. Taking $k=4$

$$D(\mathbf{A}) := \begin{array}{c} \\ x_1 \\ x_3 \\ x_2 \\ x_4 \end{array} \begin{array}{c} y_1 \quad\;\; y_3 \quad\;\; y_2 \quad\;\; y_4 \end{array} \\ \begin{bmatrix} -13/7 & 9/7 & 1/7 & 10/7 \\ 37/28 & -6/7 & -5/28 & -15/28 \\ -13/28 & 4/7 & 1/28 & 3/28 \\ 43/28 & -10/7 & 1/28 & -25/28 \end{bmatrix}$$

$B := \{x_1, x_3, x_2, x_4\}$, $N := \{y_1, y_3, y_2, y_4\}$

Now place the elements with respect to indices of variables in $B$ and $N$.

For example Here $H := \{1,2,3,4\}$, So $x_1 \in B$ and $y_1 \in N$, implies $a_{11} = -13/7$. Also $x_1 \in B$ and $y_3 \in N$ implies $a_{13} = 9/7$. Similarly placing the remaining elements we get

$$\text{Inv}(\mathbf{A}) := \begin{array}{c} \\ x_1 \\ x_2 \\ x_3 \\ x_4 \end{array} \begin{array}{c} \begin{array}{cccc} y_1 & y_2 & y_3 & y_4 \end{array} \\ \left[ \begin{array}{cccc} -13/7 & 1/7 & 9/7 & 10/7 \\ -13/28 & 1/28 & 4/7 & 3/28 \\ 37/28 & -5/28 & -6/7 & -15/28 \\ 43/28 & 1/28 & -10/7 & -25/28 \end{array} \right] \end{array}$$

**CONCLUSION**

This paper presented easy algorithms for computations of determinant and inverse of a matrix. Since the order of the given matrix has been reduced at each step while calculating its determinant, the algorithm reduces the storage requirement (as exhibited in the example). The calculation of inverses has been done using the dictionary notation which obviates the use of permutations and makes it easier to cope with in class room teaching. Ill conditioned systems can also be handled as the selection of pivots has been kept arbitrary, thus improving the numerical accuracy of the systems.

**REFERENCES**

Ahmad, F., & Khan, H. (2010). An efficient and simple Algorithm for matrix inversion. *International Journal of Technology Diffusion* , 20-27.

Augustin, C. L. (1815). Memoire sur les fonctions qui ne peuvent obtenir que deux valeurs égales et des signes contraires par suite des transpositions operées entre les variables qu'elles renferment. *de l'École polytechnique* .

Chvatal, V. (1983). *Linear Programming.* United States of America: W.H. Freeman and Company.

Cramer, G. (1750). *Introduction to the analysis of algebraic curves.*

Gauss, C. F. (1801). *Disquisitiones arithmeticae.*

Khan, H., Shah, I. A., & Ahmad, F. (2010). An efficient and generic algorithm for matrix inversion. *Internantional Jounal of Technology diffusion* , 36-41.

Maclaurin, C. (1748). *A Treatise of Algebra.* London: A. Millar and J.Nourse.

Rice, A., & Torrence, E. (2006). Lewis Caroll's Condensation method for evaluating determinants. *Math Horizons* , 12-15.